\RequirePackage{ifpdf}
\ifpdf 
\documentclass[pdftex]{sigma}
\else
\documentclass{sigma}
\fi

\input amssym.def
\newcommand{\nordbullet}{\mbox{\tiny ${\bullet\atop\bullet}$}}

\numberwithin{equation}{section}
\numberwithin{theorem}{section}
\numberwithin{proposition}{section}
\numberwithin{lemma}{section}
\numberwithin{corollary}{section}
\numberwithin{conjecture}{section}
\numberwithin{remark}{section}
\numberwithin{definition}{section}

\begin{document}


\renewcommand{\thefootnote}{$\star$}

\renewcommand{\PaperNumber}{087}

\FirstPageHeading

\ShortArticleName{The $N=1$ Triplet Vertex Operator Superalgebras:
Twisted Sector}

\ArticleName{The $\boldsymbol{N=1}$ Triplet Vertex Operator Superalgebras:\\
Twisted Sector\footnote{This paper is a
contribution to the Special Issue on Kac--Moody Algebras and Applications. The
full collection is available at
\href{http://www.emis.de/journals/SIGMA/Kac-Moody_algebras.html}{http://www.emis.de/journals/SIGMA/Kac-Moody{\_}algebras.html}}}

\Author{Dra\v{z}en ADAMOVI\'C~$^\dag$ and Antun MILAS~$^\ddag$}

\AuthorNameForHeading{D. Adamovi\'c  and A. Milas }

 \Address{$^\dag$~Department of Mathematics, University of Zagreb, Croatia}
\EmailD{\href{mailto:adamovic@math.hr}{adamovic@math.hr}}
\URLaddressD{\url{http://web.math.hr/~adamovic}}

\Address{$^\ddag$~Department of Mathematics and Statistics, University at
Albany (SUNY),\\
\hphantom{$^\ddag$}~Albany, NY 12222, USA}
\EmailD{\href{mailto:amilas@math.albany.edu}{amilas@math.albany.edu}}
\URLaddressD{\url{http://www.albany.edu/~am815139}}

\ArticleDates{Received August 31, 2008, in f\/inal form December 05,
2008; Published online December 13, 2008}

\Abstract{We classify irreducible $\sigma$-twisted modules for the
$N=1$ super triplet vertex operator  superalgebra $\mathcal{SW}(m)$
introduced recently [Adamovi\'{c}  D., Milas A.,  {\em Comm. Math. Phys.}, to appear,
\href{http://arxiv.org/abs/0712.0379}{arXiv:0712.0379}].
Irreducible
graded dimensions of $\sigma$-twisted modules are also determined.
These results, combined with our previous work in the untwisted
case, show that the $SL(2,\mathbb{Z})$-closure of the space spanned
by irreducible characters, irreducible supercharacters and
$\sigma$-twisted irreducible characters is $(9m+3)$-dimensional. We
present strong evidence that this is also the (full) space of
generalized characters for $\mathcal{SW}(m)$.  We are also able to
relate irreducible $\mathcal{SW}(m)$ characters to characters for
the triplet vertex algebra $\mathcal{W}(2m+1)$, studied in
[Adamovi\'{c} D., Milas A., {\em Adv. Math.} {\bf 217} (2008), 2664--2699,
\href{http://arxiv.org/abs/0707.1857}{arXiv:0707.1857}].}

\Keywords{vertex operator superalgebras; Ramond twisted
representations}

\Classification{17B69; 17B67; 17B68; 81R10}

\def \<{\langle}
\def \>{\rangle}
\def \t{\tau }
\def \a{\alpha }
\def \e{\epsilon }
\def \l{\lambda }
\def \L{ L}
\def \T{\mathcal T}
\def \b{\beta }
\def \om{\omega }
\def \o{\omega }
\def \c{\chi}
\def \ch{\chi}
\def \cg{\chi_g}
\def \ag{\alpha_g}
\def \ah{\alpha_h}
\def \ph{\psi_h}
\def \gi {\gamma}

\newcommand{\wta}{{\rm {wt} }  a }
\newcommand{\R}{\frak R}

\newcommand{\wtb}{{\rm {wt} }  b }

\newcommand{\g}{\frak g}

\def \ga{\gamma }
\newcommand{\hg}{\hat {\frak g} }
\newcommand{\hn}{\hat {\frak n} }
\newcommand{\h}{\frak h}
\newcommand{\wt}{{\rm {wt} }   }
\newcommand{\V}{\cal V}
\newcommand{\hh}{\hat {\frak h} }

\newcommand{\n}{\frak n}
\newcommand{\Z}{\Bbb Z}
\newcommand{\W}{\mathcal W}

\newcommand{\Zp}{{\Bbb Z}_{>0} }

\newcommand{\N}{{\Bbb Z}_{\ge 0} }
\newcommand{\C}{\Bbb C}
\newcommand{\Q}{\Bbb Q}
\newcommand{\vak}{\bf 1}

\newcommand{\la}{\langle}
\newcommand{\ra}{\rangle}

\newcommand{\NS}{\frak{ns} }
\newcommand{\triplet}{\mathcal{W}(p)}
\newcommand{\striplet}{\mathcal{SW}(m)}
\newcommand{\ssinglet}{\overline{{SM}(1)}}
\newcommand{\hf}{\mbox{$\frac{1}{2}$}}
\newcommand{\thf}{\mbox{$\frac{3}{2}$}}

\section{Introduction}

Many constructions and results in vertex algebra theory are easily extendable to the setup of vertex superalgebras by simply adding adjective ``super''. Still, there are results that deviate from this ``super-principle'' and new ideas are needed compared to the non-super case. For example, modular invariance for vertex operator superalgebras requires inclusion of supercharacters of (untwisted) modules, and more importantly, the characters of $\sigma$-twisted modules \cite{FFR}, where $\sigma$ is the canonical parity automorphism. Since the construction and classif\/ication of $\sigma$-twisted modules is more or less independent of the untwisted construction, many aspects of the theory need to be reworked for the twisted modules (e.g., twisted Zhu's algebra~\cite{Xu}). In fact, even for the free fermion vertex operator superalgebra, construction of $\sigma$-twisted modules is far from being trivial (see~\cite{FFR} for details).

Present work is a natural continuation of our very recent paper
\cite{AdM-str}, where we introduced a~new family of $C_2$-cof\/inite
$N=1$ vertex operator superalgebra that we call the {\em
supertriplet} family~$\mathcal{SW}(m)$, $m \geq 1$. In this
installment we focus on $\sigma$-twisted $\mathcal{SW}(m)$-modules
and their irreducible characters. The $\sigma$-twisted $\striplet$-modules are usually called modules in the {\em Ramond sector}, while untwisted (ordinary) modules are referred to as modules in the {\em Neveu--Schwarz sector}, in parallel with two distinct $N=1$ superconformal algebras.

Of course, Ramond sector has been studied for various (mostly
rational) vertex operators superalgebras (e.g.,
\cite{Li-tw,D2,M}, etc.). We should also point out
that Ramond twisted representation of certain vertex
$\mathcal{W}$-algebras were recently investigated in  \cite{Ar-t}
and \cite{KW-t}.

Here are our main results:

\begin{theorem}
The twisted Zhu's algebra $A_{\sigma}(\striplet)$ is finite-dimensional with precisely $2m\!{+}1$, non-isomorphic,  $\mathbb{Z}_2$-graded irreducible modules.
 Consequently, $\striplet$ has $2m+1$, non-iso\-mor\-phic, $\mathbb{Z}_2$-graded irreducible $\sigma$-twisted modules.
\end{theorem}

By using embedding structure of $N=1$ Feigin--Fuchs modules \cite{IK3} we can also easily compute the characters of $\sigma$-twisted $\striplet$-modules. Equipped with these formulas, results from \cite{AdM-str} about untwisted modules and supercharacters, and transformation formulas for classical Jacobi theta functions we are able to prove:
\begin{theorem}
The $SL(2,\mathbb{Z})$-closure of the space of characters, (untwisted) supercharacters and $\sigma$-twisted $\striplet$ characters is $(9m+3)$-dimensional.
\end{theorem}

Conjecturally, this is also the space of certain {\em generalized characters} studied in \cite{Miy}  (strictly speaking pseudotraces are yet-to-be def\/ined in the setup of $\sigma$-twisted modules). A closely
related conjecture is

\begin{conjecture} Let $Z(A)$ denote the center of associative algebra $A$ and $T(A)=A/[A,A]$ the trace group of $A$. Then
\[
\dim (Z(A(\striplet)))+\dim (Z(A_{\sigma}(\striplet)))=9m+3,
\]
and
\[
\dim (T(A(\striplet)))+\dim (T(A_{\sigma}(\striplet)))=9m+3,
\]
where $A(\striplet)$ is the untwisted Zhu's associative algebra of
$\striplet$ (cf.~{\rm \cite{AdM-str}}).
\end{conjecture}
The conjecture would follow if we knew more about structure of logarithmic $\striplet$-modules.

Here is a short outline of the paper.  In Section~\ref{prelim} we
recall the standard results about the $\sigma$-twisted Zhu's
algebra, the free fermion vertex superalgebra $F$, and the
construction of the $\sigma$-twisted $F$-module(s). In Section
\ref{singlet-ramond} we focus on $\sigma$-twisted modules for the
super singlet vertex algebra introduced also in \cite{AdM-str}.
Sections \ref{free-fields} and \ref{classification} deal with the
construction and classif\/ication of $\sigma$-twisted
$\striplet$-modules. Finally, in Section \ref{characters} we derive
characters formulas for irreducible $\sigma$-twisted
$\striplet$-modules, discuss modular invariance, and relate our
characters with the cha\-rac\-ters for the triplet vertex algebra
$\mathcal{W}(2m+1)$ (see for instance \cite{AdM-triplet} and \cite{FGST}).

Needless to say, this paper is largely continuation of \cite{AdM-triplet} and especially \cite{AdM-str}.  Thus the reader is strongly encouraged to consult \cite{AdM-str}, where we studied $\striplet$ in great details.

\section{Preliminaries}
\label{prelim} Let $V=V_0 \oplus V_1$ be a vertex operator
superalgebra, where as usual $V_0$ is the even and $V_1$ is the odd
subspace (cf.~\cite{FFR,HM,K}). Every vertex operator superalgebra has the canonical
parity automorphism $\sigma$, where $\sigma_{V_0}=1$ and
$\sigma_{V_1}=-1$. This leads to the notion of $\sigma$-twisted
$V$-modules, well recorded in the literature (see for example
\cite{FRW} and \cite{Xu}). As in the untwisted case, a large part of
representation theory of $\sigma$-twisted $V$-modules can be
analyzed via the $\sigma$-twisted Zhu's algebra whose construction
we recall here.

Within the same setup, consider the subspace $O(V)
\subset V$, spanned by elements of the form
\[
{\rm Res}_x \frac{(1+x)^{{\rm deg}(u)}}{x^2} Y(u,x)v,
\]
where $u \in V$ is homogeneous. It can be easily shown that
\[
{\rm Res}_x \frac{(1+x)^{{\rm deg}(u)}}{x^n} Y(u,x)v \in O(V) \qquad \mbox{for} \quad n \geq 2.
\]
Then, the vector space $A_\sigma(V)=V/O(V)$ is equipped with
an associative algebra structure via
\[
u * v={\rm Res}_x \frac{(1+x)^{\rm deg(u)}}{x} Y(u,x)v
\]
(see~\cite{Xu}). An important dif\/ference between the untwisted
associative algebra $A(V)$ and $A_\sigma(V)$ is that the latter is
$\mathbb{Z}_2$-graded, so
\[
A_\sigma(V)=A^0_\sigma(V) \oplus A^1_\sigma(V).
\]
It is also not clear that $A_{\sigma}(V) \neq 0$, while $A(V)$ is
always nontrivial. We shall often use $[a] \in A_{\sigma}(V)$ for
the image of $a \in V$ under the natural map $V \longrightarrow
A_{\sigma}(V)$.

The result we need is \cite{Xu}:

\begin{theorem} \label{zhu-corr} There is a one-to-one correspondence between irreducible
$\mathbb{Z}_{\geq 0}$-gradable $\sigma$-twisted $V$-modules and
irreducible $A_{\sigma}(V)$-modules.
\end{theorem}
In the theorem there is no reference to graded $A_{\sigma}(V)$-modules.
For practical purposes we shall need a slightly dif\/ferent version of the above theorem, because some modules are more natural if considered as $\mathbb{Z}_2$-graded modules (shorthand,  {\em graded} modules).

\begin{theorem} \label{zhu-graded-corr} There is a one-to-one correspondence between graded irreducible $\mathbb{Z}_{\geq 0}$-gradable $\sigma$-twisted $V$-module
and graded irreducible $A_{\sigma}(V)$-module.
\end{theorem}

\begin{proof}
The proof mimics the non-graded case, so we shall omit
details. As in the non-graded case,  by applying Zhu's theory, from
an irreducible graded $A_{\sigma}(V)$-module $U$ we construct a~$\sigma$-twisted graded $V$-module $L(U)$ (but of course in the
process of getting $L(U)$ will be moding out by the maximal {\em
graded} submodule). On the other hand,  if $M$ is an irreducible
graded $V$-module, then the top component $\Omega(M)$ is clearly a
graded $A_{\sigma}(V)$-module. Suppose that $\Omega(M)$ is not
graded irreducible, therefore there is a graded submodule
$\Omega(M')$. Then  we let $M'=\mathcal{U}(V[\sigma])\Omega'(M) $,
where $\mathcal{U}(V[\sigma])$ is the enveloping algebra of the Lie
algebra $V[\sigma]$.
 But $M'$ is a proper graded submodule of $M$, so we get  a contradiction.
\end{proof}

Our goals are to describe the structure of
$A_\sigma(\overline{SM(1)})$, where $\overline{SM(1)}$ is the super
singlet vertex algebra \cite{AdM-str}, and to discuss
$A_\sigma(\striplet)$ (in fact, we have a very precise conjecture
about the structure of $A_\sigma(\striplet)$). Let us recall (see
\cite{AdM-str} for details) that both $\ssinglet$ and $\striplet$
are $N=1$ superconformal vertex operator superalgebras, with the
superconformal vector $\tau$. In other words, if we let
$Y(\tau,x)=\sum\limits_{n \in \mathbb{Z}+1/2} G(n)x^{-n-3/2}$, then $G(n)$
and $L(m)$ close the $N=1$ Neveu--Schwarz superalgebra.

 It is not
hard to see that the following hold:{\samepage
\begin{gather}
 (L(-m-2)+2L(-m-1)+L(-m))v \in O(V), \qquad m \geq 2, \nonumber\\
(L(0)+L(-1))v \in O(V), \nonumber\\
 L(-m)v \equiv (-1)^m((m-1)(L(-2)+L(-1))+L(0))v \mod O(V), \qquad m \geq 2,\nonumber\\
   \sum_{n \geq 0} {\frac{3}{2} \choose n} G(-3/2-i+n) v \in O(V),
\qquad i \geq 1,   \label{ramond-relations}
 \end{gather}
  where in all formulas $v$ is a vector in $N=1$
vertex operator superalgebra $V$.}

To illustrate how to use twisted Zhu's algebra, let us classify
$\mathbb{Z}$-graded $\sigma$-twisted modules for the (neutral) free
fermion vertex operator superalgebra $F$, used in \cite{AdM-str}.
The next result is known so we will not provide its proof here.

 \begin{proposition}   We have $A_{\sigma}(F) \cong {\C}[x]/\big(x^2-\frac{1}{2}\big)$ where $x=[\phi(-1/2){\bf 1}]$,
an odd generator in the associative algebra.
\end{proposition}

 This result in particular implies there are precisely two
irreducible $\sigma$-twisted $F$-modules: $M^{\pm}$. These two
modules can be constructed explicitly. As vector spaces
\[
M^{\pm}=
\oplus_{n \geq 0} \Lambda^n (\phi({-1}),\phi({-2}),\dots),
\]
 where
$\Lambda^*$ is the exterior algebra, which is also a
$\mathbb{Z}_{\geq 0}$-graded module for the Clif\/ford algebra $K$
spanned by $\phi({n})$, $n \in \mathbb{Z}$ with anti-bracket
relations $\{ \phi(m),\phi(n) \}=\delta_{m+n,0}$. The only
dif\/ference between $M^+$ and $M^-$  is in the action of $\phi(0)$ on
the one-dimensional top subspace $M^\pm(0) $. More precisely we have
$\phi(0)|_{M^\pm(0)}=\pm \frac{1}{\sqrt{2}}$. However, notice that
$M^{\pm}$ are not $\mathbb{Z}_2$-graded thus it is more natural to examine graded
$A_{\sigma}(F)$-modules. There is a unique such module (up to parity
switch), spanned  by ${\bf 1}^{\R}$ and $\phi(0){\bf 1}^{\R}$.
 Thus
\[
M=\oplus_{n \geq 0} \Lambda^n(\phi(0), \phi({-1}),\phi({-2}),\dots) =
M^+ \oplus M^-.
\]
Then   ${\bf 1}^{\R}$ is a cyclic vector in $M$, i.e., $M= K. {\bf
1}^{\R}$. Moreover, $ M^{\pm} = K. {\bf 1}^{\pm}$,
where
\[
 {\bf 1}^{\pm} = {\bf 1} ^{\R} \pm \sqrt{2} \phi(0) {\bf 1} ^{\R}.
\]

 Next we describe the  twisted vertex operators
 \[
 \overline{Y} :  \ \ F
\otimes M \longrightarrow M[[x^{1/2},x^{-1/2}]] .
\]
 Details are
spelled out in \cite{FFR}, here we only give the explicit formula.
Def\/ine f\/irst
\[
Y\big(\phi\big({-}n-\tfrac{1}{2}\big){\bf 1},x\big)=\frac{1}{n!}\left(\frac{d}{dx}\right)^{n}\left(\sum_{m
\in \mathbb{Z}} \phi(m)x^{-m-1/2}\right),
\] acting on $M$ and use
(fermionic!) normal ordering ${\nordbullet } \ \ {\nordbullet }$ to
def\/ine
\[
Y\big(\phi\big({-}n_1-\tfrac{1}{2}\big)\cdots \phi\big({-}n_k-\tfrac{1}{2}\big){\bf 1},x\big)=\nordbullet Y\big(\phi\big({-}n_1-\tfrac{1}{2}\big),x\big) \cdots Y\big(\phi\big({-}n_k-\tfrac{1}{2}\big),x\big) \nordbullet,
\]
and extend $Y$ by linearity on all of $F$ (see \cite{FFR} and \cite{FRW} for details, especially about normal ordering).
Then, we let
\[
\overline{Y}(v,x):=Y(e^{\Delta_x}v,x),
\]
where
\begin{gather*}
  \Delta_x=\frac{1}{2} \sum_{m,n \geq 0}  C_{m,n} \phi\big(m+\tfrac{1}{2}\big)\phi\big(n+\tfrac{1}{2}\big) x^{-m-n-1}, \nonumber \\
  C_{m,n}=\frac{1}{2}\frac{m-n}{m+n+1}{-1/2 \choose m}{-1/2 \choose n}.\nonumber
\end{gather*}
Let us f\/ix the Virasoro generator
\[
\omega_s=\tfrac{1}{2} \phi\big({-}\tfrac{3}{2}\big)\phi\big({-}\tfrac{1}{2}\big){\bf 1}.
\]
Then $e^{\Delta_{x}}\omega_s=\omega_s+\frac{1}{16}x^{-2}{\bf 1}$.

The following lemma will be important in the rest of the paper.

\begin{lemma} \label{cmn} Let
\[
G(x_1,x_2)=\sum_{m=0}^\infty \sum_{n=0}^\infty C_{m,n} x_1^m x_2^n,
\]
where $C_{m,n}$ are as above. Then
\[
G(x_1,x_2)=\frac{1}{2}\left(\frac{(1+x_1)^{1/2}(1+x_2)^{-1/2}-1}{x_2-x_1}+\frac{(1+x_1)^{-1/2}(1+x_2)^{1/2}-1}{x_2-x_1}\right)  \in \mathbb{Q}[[x_1,x_2]].
\]
\end{lemma}

\begin{proof} The lemma follows directly from the identity
\[
(x_2-x_1)G(x_1,x_2)=\tfrac{1}{2}(1+x_1)^{1/2}(1+x_2)^{-1/2}+\tfrac{1}{2}(1+x_1)^{-1/2}(1+x_2)^{1/2}-1.
\]
which can be easily checked by expanding both sides as power
series in $x_1$ and $x_2$.
 \end{proof}

\section{Highest weight representations of the Ramond algebra}

 The $N=1$ Ramond  algebra $\R$ is the
inf\/inite-dimensional Lie superalgebra
\[
{\R} = \bigoplus_{n \in {\Z} } {\C} L(n) \bigoplus \bigoplus_{m
\in  {\Z} }  {\C}G(m) \bigoplus {\C}C
\]
 with commutation relations ($m,n \in {\Z} $):
\begin{gather*}
[L(m),L(n)] =(m-n)L({m+n})+\delta_{m+n,0}\frac{m^3-m}{12}C,
\nonumber  \\   
[G(m),L(n)]   =
\left(m-\frac{n}{2}\right)G({m+n}),    \\  
\{ G(m ),G(n) \} =
2 L(m+n)+\frac{1}{3}(m^2 -1/4)\delta_{m+n,0}C,
\\   [L(m),C] = 0, \qquad [G(m),C] = 0.
\nonumber
\end{gather*}

The representation theory of the $N=1$ Ramond algebra has been intensively
studied f\/irst in \cite{KW-85} and other papers (cf.~\cite{MR,Dor}, etc.).

Assume that $(c,h) \in {\C}^2$ such that $  24 h \ne   c $. Let
$L^{\R} (c,h)^{\pm}$ denote the irreducible  highest weight
$\R$-module generated by the highest weight vector $v_{c,h}
^{\pm}$ such that
\begin{gather*}
  G(n) v_{c,h} ^{\pm} = \pm \sqrt{h -c /24} \,\delta_{n,0} v_{c,h}
^{\pm}, \qquad
 L(n) v_{c,h} ^{\pm} = h \delta_{n,0} v_{c,h} ^{\pm} \qquad (n \ge
0). \end{gather*}
These modules can be considered as irreducible $\sigma$-twisted
modules for the Neveu--Schwarz vertex operator superalgebra (cf.~\cite{Li-tw}). Since $L^{\R} (c,h)^{\pm}$ are not ${\Z}_2$-graded
(notice that $v_{c,h}^\pm$ are eigenvectors for $G(0)$), it is more
useful to consider  graded modules. It is not hard to show that the
direct sum
\[
  L^{\R} (c,h) = L^{\R} (c,h)^{+} \oplus L^{\R} (c,h)^{-}
\]
is in fact a ${\Z}_2$-graded $\R$-module. Indeed, for
$\mathbb{Z}_2$-graded subspaces take
\[
L^{\R}(c,h)_0=U(\R)\left(v^{+}_{c,h}+\frac{1}{\sqrt{h-c/24}}v^{-}_{c,h}\right)
\]
and
\[
 L^{\R}(c,h)_1=U(\R)\left(v^{+}_{c,h}-\frac{1}{\sqrt{h-c/24}}v^{-}_{c,h}\right)
\]
 Since
$L^{\R} (c,h)$ does not contain non-trivial ${\Z}_2$-graded
submodules, we shall say that this module is ${\Z}_2$-graded
irreducible.

\begin{remark}
Details about construction of  $L^{\R} (c,h)$ and its relation to
irreducible non-graded modules can be found in~\cite{Dor} and~\cite{IK2}. In Section~\ref{free-fields} we shall present   free
f\/ields realization of modules $L^{\R} (c,h)$ and $L^{\R}
(c,h)^{\pm}$.
\end{remark}

For  $i, n, m \in {\N}$, let
\begin{gather*}
c_{2m+1,1} = \frac{3}{2} - \frac{12 m ^2 }{2m +1}, \qquad h^{1,3} = 2m +\frac{1}{2}, \\
 h^{2 i +2,2n+1}:= \frac{ (2 i + 2
 - (2n+1) (2m+1) ) ^2 - 4m^2 }{8 (2m
+1)} + \frac{1}{16}.
\end{gather*}

Let $ L(c_{2m+1,1},h)$ be the irreducible highest weight module for
the Neveu--Schwarz algebra with central charge $c$  and highest
weight~$h$. Then $L(c_{2m+1,1},0)$ is a simple vertex operator
superalgebra, $ L^{\R} (c,h)^{\pm} $ are irreducible
$\sigma$-twisted  $L(c_{2m+1,1},0)$-modules, and $ L^{\R} (c,h) $
is ${\Z}_2$-graded irreducible $\sigma$-twisted
$L(c_{2m+1,1},0)$-module.

\subsection{Intertwining operators among twisted modules}

If $V$ is a vertex operator superalgebra and $W_1$, $W_2$ and $W_3$ are three $V$-modules then we consider the space of
intertwining operators ${ W_3 \choose W_1 \ W_2}$. It is perhaps less standard to study intertwining operators between twisted
$V$-modules, so we recall the def\/inition here  (see \cite[p.~120]{Xu}).

\begin{definition} Let $W_1$, $W_2$ and $W_3$ be $\sigma_i$-twisted
$V$-modules, respectively, where $\sigma_i$ is a f\/inite order automorphism
of order $\nu_i$, with common period $T$.  An intertwining operator of type
${W_3 \choose W_1 \ \ W_2 }$ is a linear map
\[
\mathcal{Y}( \cdot, x) : \ \ W_1 \longrightarrow  {\rm End}(W_2,W_3)\{ x \},
\]
such that
\[
\mathcal{Y}(L(-1)w,x)=\frac{d}{dx}\mathcal{Y}(w,x), \qquad w \in W_1
\]
and Jacobi identity holds 
\begin{gather*} \frac{1}{T x_0} \sum_{p=0}^{T-1}
\delta \left(\left(\frac{x_1-x_2}{x_0}  \right)^{-1/T} \omega_T^p
\right){Y}^{\nu_3}(\sigma_1^p v,x_1)
\mathcal{Y}(w_1,x_2) \nonumber \\
\qquad{} -(-1)^{ij}
\frac{1}{T x_0} \sum_{p=0}^{T-1}\delta \left(\left(\frac{-x_2+x_1}{x_0}
\right)^{-1/T} \omega_T^p \right) \mathcal{Y}(w_1,x_2){Y}^{\nu_2}(\sigma_1^p
v,x_1) \nonumber \\
\qquad{}=\frac{1}{T x_2} \sum_{p=0}^{T-1} \delta
\left(\left(\frac{x_1-x_0}{x_2}  \right)^{1/T} \omega_T^p
\right)\mathcal{Y}(Y^{\nu_1}(\sigma_2^p v,x_0)w_1,x_2), 
\end{gather*}
 where $v \in V_i$, $w_1 \in (W_1)_j$, $i,j \in \mathbb{Z}_2$
and $\omega_T$ is a primitive $T$-th root of unity.
\end{definition}

We shall also need the following result on the fusion rules. The
proof is completely analogous to that of Proposition~4.1 of~\cite{AdM-str} (see also~\cite{IK1}).

\begin{proposition} \label{fusion-R}
For every $i=0,\dots,m-1$ and $n \geq 1$ we have: the space
\[
I   {L^{\R}(c_{2m+1,1},h)^{\pm} \choose L(c_{2m+1,1},h^{1,3}) \
L^{\R}(c_{2m+1,1},h^{2i+2,2n+1})^{\pm}}
\]
 is nontrivial only if $h
\in \{h^{2i+2,2n-1}, h^{2i+2,2n+1},h^{2i+2,2n+3} \}$, and
\[
I
{L^{\R}(c_{2m+1,1},h)^{\pm}  \choose  L(c_{2m+1,1},h^{1,3}) \
L^{\R}(c_{2m+1,1},h^{2i+2,1}) ^{\pm}}
\] is nontrivial only if
$h=h^{2i+2,3}.$

Similarly, for every $i=0,\dots,m-1$ and $n \geq 2$ we have: the space
\[
I  {L^{\R}(c_{2m+1,1},h)^{\pm} \choose L(c_{2m+1,1},h^{1,3}) \
L^{\R}(c_{2m+1,1},h^{2i+2,-2n+1})^{\pm}}
\]
 is nontrivial only if $h
\in \{ h^{2i+2,-2n-1},h^{2i+2,-2n+1}, h^{2i+2,-2n+3} \}$, and
\[
I { L^{\R}(c_{2m+1,1},h)^{\pm} \choose L(c_{2m+1,1},h^{1,3}) \
L^{\R}(c_{2m+1,1},h^{2i+2,-1})^{\pm} }
\] is nontrivial only if $h
\in \{h^{2i+2,-3},h^{2i+2,-1} \}$.

We have analogous result for fusion rules of type
\[
{L^{\R}(c_{2m+1,1},h)  \choose L(c_{2m+1,1},h^{1,3}) \
L^{\R}(c_{2m+1,1},h^{2i+2,2n+1})  }.
\]
\end{proposition}

 \begin{proof} The proof goes along the lines in~\cite{AdM-str} (cf.~\cite{M0}), with some minor modif\/ications due to
twisting. In fact, in order to avoid the  twisted version of
Frenkel--Zhu's formula we can proceed in a more straightforward
fashion. Because all modules in question are irreducible and because
all intertwining operators we are interested in are of the form
${W_3 \choose W_1 \ \ W_2}$, where $W_1$ is untwisted (so
$\nu_1=0$), the left-hand side in the Jacobi identity looks like the
ordinary Jacobi identity for intertwining operators. Thus we can use
commutator formula and the null vector conditions
\begin{gather*}
\big({-}L(-1)G\big({-}\tfrac{1}{2}\big)+(2m+1)G\big({-}\tfrac{3}{2}\big)\big)v_{1,3}=0,\\
G\big({-}\tfrac{1}{2}\big)\big({-}L(-1)G\big({-}\tfrac{1}{2}\big)+(2m+1)G\big({-}\tfrac{3}{2}\big)\big)v_{1,3}=0,
\end{gather*}
which hold in $L(c_{2m+1,1},h^{1,3})$ (here $v_{1,3}$ is the highest weight vector),
to study the matrix coef\/f\/icient
\[
f(x)=\langle w',\mathcal{Y}(v_{1,3},x) w \rangle,
\]
where $w$ and $w'$ are highest weight vectors in appropriate
modules. This leads to dif\/ferential equations for $f(x)$, which can
be solved. The general solution is a linear combination of power
functions $x^s$, where $s$ is a rational number. The rest follows by
interpreting $s$  in terms of conformal weights for the three
modules involved in~$\mathcal{Y}$.
\end{proof}

\section[$\sigma$-twisted modules for the super singlet algebra $\overline{SM(1)}$]{$\boldsymbol{\sigma}$-twisted modules for the super singlet algebra $\boldsymbol{\overline{SM(1)}}$}
\label{singlet-ramond}

In this section  $\sigma_2$ will  denote the parity automorphism of~$F$. Recall also the Heisenberg vertex operator algebra $M(1)$.
Then, as in \cite{AdM-str} we equip the space $M(1) \otimes F$
with a vertex superalgebra structure such that the total central
charge is $\frac{3}{2}(1-\frac{8m^2}{2m+1})$. We def\/ine the following superconformal and conformal
vectors:
\begin{gather*}
\tau = \frac{1}{ \sqrt{2 m +1 }} \left( \alpha(-1) {\vak}
\otimes \phi \big({-}\hf\big) {\vak} + 2 m {\vak} \otimes  \phi \big({-\thf}\big)
{\vak}
\right),  \\
\omega = \frac{1}{2 (2m+1)} ( \alpha(-1) ^{2} +  2 m \alpha(-2) )
{\vak} \otimes {\vak} +    {\vak} \otimes  {\omega}^{(s)}.
\end{gather*}

Recall also the singlet superalgebra $\overline{SM(1)}$ obtained
as the kernel of the screening operator
\[
\widetilde{Q}={\rm Res}_x Y\big(e^{\frac{-\alpha}{2m+1}} \otimes \phi\big({-}\tfrac{1}{2}\big),x\big)
\] acting from $M(1) \otimes F$ to
$M(1) \otimes e^{-\frac{\alpha}{2m+1}} \otimes F$, where $ \langle
\alpha,\alpha \rangle=2m+1$. The vertex operator superalgebra
$\overline{SM(1)}$ is generated by $\tau=G\big({-}\tfrac{3}{2}\big){\bf 1}$ and
$H=Qe^{-\alpha}$, where
\[
Q={\rm Res}_{x} Y\big(e^{\alpha} \otimes \phi\big({-}\tfrac{1}{2}\big),x\big).
\]
We will also use $\omega$ and $\widehat{H}$, where the latter is
proportional to $G\big({-}\tfrac{1}{2}\big)H$ (for details see \cite{AdM-str}).

Consider the automoprhism $\sigma= 1 \otimes \sigma_2$ of $M(1) \otimes F$, acting
nontrivially on the second tensor factor. This automorphism
plainly preserves $\overline{SM(1)}$, thus we can study
$\sigma$-twisted $\overline{SM(1)}$-modules. It is clear that
$M(1,\lambda) \otimes M$ is a graded $\sigma$-twisted
$\overline{SM(1)}$-module, and $M(1,\lambda) \otimes M^{\pm}$ are
$\sigma$-twisted $\overline{SM(1)}$-modules, where $M(1,\lambda)$ is as in
\cite{AdM-triplet} and  $\lambda \in \mathbb{C}$.

We would like to classify irreducible $\overline{SM(1)}$-modules by
virtue of Zhu's algebra. As usual, we denote by $[a] \in
A_{\sigma}(\overline{SM(1)})$ the image of $a \in \overline{SM(1)}$.

\begin{lemma} \label{formule}
Let $v_{\lambda}^{\pm} = v_{\lambda}\otimes {\bf 1}^{\pm}$ be the
highest weight vector in $M(1,\lambda) \otimes M^{\pm}$, where
$\langle \alpha,\lambda \rangle=t$. Then the twisted generators
act as
\begin{gather*}
G(0) \cdot v_{\lambda}^{\pm}=\pm \frac{t-m}{\sqrt{2(2m+1)}} v_{\lambda} ^{\pm},\\
L(0) \cdot v_{\lambda}^{\pm}=\left(\frac{t(t-2m)}{2(2m+1)}+\frac{1}{16} \right) v_{\lambda} ^{\pm},\\
 H(0) \cdot v_{\lambda}^{\pm}=\pm \frac{1}{\sqrt{2}}{t-1/2 \choose 2m} v_{\lambda} ^{\pm},\\
 \widehat{H}(0) \cdot v_{\lambda} ^{\pm} = \frac{m-t}{2m+1} {t-1/2 \choose 2m} v_{\lambda} ^{\pm}.
 \end{gather*}
\end{lemma}

\noindent \begin{proof} Recall from \cite{AdM-str} the formulas
\begin{gather*}
H=S_{2m}(\alpha) \otimes \phi\big({-}\tfrac{1}{2}\big) +S_{2m-1} \otimes \phi\big({-}\tfrac{3}{2}\big) +  \cdots + 1 \otimes \phi\big({-}2m-\tfrac{1}{2}\big),\\
 \widehat{H}=\phi(1/2)S_{2m+1}(\alpha)\phi\big({-}\tfrac{1}{2}\big) + w \\
 \hphantom{\widehat{H}}{} =S_{2m+1}(\alpha) + S_{2m -1} (\a) \otimes \phi\big({-}\tfrac{3}{2}\big) \phi\big({-}\tfrac{1}{2}\big) + \cdots + {\bf 1} \otimes \phi\big({-}2m -\tfrac{1}{2}\big) \phi\big({-}\tfrac{1}{2}\big).
 \end{gather*}
As in \cite{AdM-str} we have
\[
S_{r}(\alpha)(0) \cdot v_{\lambda} ^{\pm}={t \choose r} v_{\lambda} ^{\pm}.
\]
Then we get
\[
H(0) \cdot v_{\lambda} ^{\pm}=\pm \frac{1}{\sqrt{2}}\sum_{n=0}^{2m} {t \choose n}{-1/2 \choose 2m-n} v_{\lambda} ^{\pm}={\pm} \frac{1}{\sqrt{2}}{t-1/2 \choose 2m}
v_{\lambda} ^{\pm}.
\] The last formula is proven similarly by using
$e^{\Delta_x}$ operator. \end{proof}

Here are the main result of this section.
\begin{theorem} \label{relacije-sin} The Zhu's algebra $A_{\sigma}(\overline{SM(1)})$ is an associative algebra generated by
 $[\tau]$, $[\omega]$, $[H]$ and $[\widehat{H}]$, where the following relations hold (here $[\omega]$ is central):
\begin{gather}
 [\tau] ^2 = [\omega] - \frac{c_{2m+1,1}}{24},\nonumber\\
 [\tau]*[H] =[H]*[\tau]=\frac{-\sqrt{2m+1}}{2} [\widehat{H}],\nonumber\\
 [H] * [H] = \frac{1}{2} { \sqrt{2(2m+1)} [\tau] + {m-1/2} \choose 2
m}  ^2,  \nonumber\\
\label{important}
[H]*[H]=C_m \prod_{i=0}^{m-1}\big([\omega]-h^{2i+2,1}\big)^2,
\qquad C_m = \frac{2 ^{2m-1} (2m+1) ^{2m}}{ (2m)!^2},\\
 [\widehat{H}]*[\widehat{H}]=\frac{4}{2m+1} C_m \left([\omega] - \frac{c_{2 m+1,1}}{24}\right) \prod_{i=0}^{m-1}\big([\omega]-h^{2i+2,1}\big)^2. \nonumber
\end{gather}
In particular, Zhu's algebra is  commutative.
\end{theorem}

\begin{proof} First notice (cf. \cite{AdM-str}) that
\[
G\big({-}\tfrac{1}{2}\big){F}=-\sqrt{2m+1} \phi\big({-}\tfrac{1}{2}\big)F=-\sqrt{2m+1} \widehat{F}
\] and consequently
\[
G\big({-}\tfrac{1}{2}\big)H=-\sqrt{2m+1} \widehat{H}.
\]
We also have the relation
\[
\big((2m+1)G\big({-}\tfrac{3}{2}\big)-L(-1)G\big({-}\tfrac{1}{2}\big)\big)H=0,
\]
because $H=v_{1,3}$ is a highest weight vector.
By using (\ref{ramond-relations}) we obtain
\begin{gather*}
[\tau]*[H]=\big[G\big({-}\tfrac{3}{2}\big)H\big]+\tfrac{3}{2}\big[G\big({-}\tfrac{1}{2}\big)H\big]
=\left[\frac{1}{2m+1}L(-1)G\big({-}\tfrac{1}{2}\big)H\right]+\tfrac{3}{2}\big[G\big({-}\tfrac{1}{2}\big)H\big]\\
\phantom{[\tau]*[H]}{}
=\tfrac{1}{2}\big[G\big({-}\tfrac{1}{2}\big)H\big]=-\frac{\sqrt{2m+1}}{2}[\widehat{H}].
\end{gather*}
On the other hand, by using skew-symmetry we also have
\begin{gather*}
 [H]*[\tau]=-\big[{\rm Res}_x x^{-1} e^{x L(-1)}Y(\tau,-x)(1+x)^{2m+\frac{1}{2}}H\big] \nonumber \\
\hphantom{[H]*[\tau]}{} =-\big[{\rm Res}_x x^{-1} e^{-x L(0)}Y(\tau,-x)(1+x)^{2m+\frac{1}{2}}H\big] \nonumber \\
\hphantom{[H]*[\tau]}{}=-\big(\big[G\big({-}\tfrac{3}{2}\big)H\big]-\big[\big(2m+\tfrac{1}{2}-(2m+1)\big)
G\big({-}\tfrac{1}{2}\big)H]\big) \nonumber \\
\hphantom{[H]*[\tau]}{}=-\big[G\big({-}\tfrac{3}{2}\big)H\big]-\tfrac{1}{2}\big[G\big({-}\tfrac{1}{2}\big)H\big]=
\tfrac{1}{2}\big[G\big({-}\tfrac{1}{2}\big)H\big].
\end{gather*} Combined, we obtain
\begin{gather*}
[\tau]*[H]=[H]*[\tau] .\tag*{\qed}
\end{gather*}\renewcommand{\qed}{}
\end{proof}

Notice that relation (\ref{important}) can be written as
\[
 [H] * [H] = C_m \prod_{i=0}^{m-1}\left([\tau]^2-\frac{(2i+1-2m)^2}{8(2m+1)} \right)^2,  \qquad C_m \neq 0.
 \]

Let ${\C}[a,b]$ denote the ${\Z}_2$-graded complex commutative
associative algebra generated by odd vectors $a$, $b$.

\begin{theorem} \label{structura-singlet}
The associative algebra $A_{\sigma}(\overline{SM(1)})$ is isomorphic to the
${\Z}_2$-graded commutative associative algebra
\[
{\C}[a,b] / \la
H(a,b) \ra,
\]
 where $\la H(a,b) \ra$ is (two-sided) ideal in
${\C}[a,b]$, generated by
\[
H(a,b) = b^2 - C_m \prod_{i=0}^{m-1}\left(a^2-\frac{(2i+1-2m)^2}{8(2m+1)}\right)^2.
\]
\end{theorem}
\begin{proof} The proof is similar to that of Theorem 6.1 of
\cite{A-2003}. First we notice that we have a~surjective
homomorphism
\begin{gather*}
 \Phi : \ \ {\C}[a,b]   \rightarrow  A_{\sigma}(\overline{SM(1)}), \qquad
 a \mapsto   [\tau],\qquad
b \mapsto   [H]  .
\end{gather*}

It is easy to see that $\mbox{Ker} \,  \Phi$ is a ${\Z}_2$-graded
ideal. We shall now prove that $\mbox{Ker} \,  \Phi =\la H(a,b) \ra$.

Evidently the generating element $H(a,b)$ is even, so Theorem~\ref{relacije-sin} gives $\la H(a,b) \ra \subset \mbox{Ker} \, \Phi$.

Assume now that $K(a,b) \in \mbox{Ker} \,  \Phi$.
By using division algorithm   we get
\[
K(a,b) = A(a,b) H(a,b) + R(a,b),
\]
where
$ A(a,b), R(a,b) \in {\C}[a,b]$
and  $R(a,b)$  has degree at most $1$ in $b$. Assume that $R(a,b)
\ne 0$. Then
$R(a,b) = A(a) b + B(a)$ for certain polynomials $A,  B \in
{\C}[x]$.  We also notice that $R(a,b) \in \mbox{Ker} \, \Phi$. Since
$\mbox{Ker} \, \Phi$ is ${\Z}_2$-graded ideal, we can assume that
$R(a,b)$ is homogeneous. If  $R(a,b)$ is an even element we have
that $A$ has odd degree and $B$ has even degree, and therefore
\begin{gather}
 \deg (B) - \deg(A) \qquad \mbox{is an odd natural number.}
\label{odd-r}
\end{gather}
The case when $R(a,b)$ is odd element again leads to formula~(\ref{odd-r}).

 As in \cite{A-2003} we now shall evaluate $R(a,b)$ on
$A_{\sigma}(\overline{SM(1)})$-modules and get
\[
   A \left( \frac{t-m}{\sqrt{2(2m+1)}} \right) {t-1/2 \choose
2m} +    \sqrt{2} B \left(\frac{t-m}{\sqrt{2(2m+1)}}\right) = 0
\qquad \forall \ t \in {\C}.
\]
This implies that $\deg(B) - \deg(A) = 2m$. Contradiction.
  Therefore $R(a,b) = 0$ and $K(a,b) \in \la H(a,b)
\ra$. The proof follows. \end{proof}

\begin{remark}
By using the same arguments as  in  \cite{A-2003} and
\cite{AdM-str}, we can conclude that every irre\-ducible
$\mathbb{Z}_{\geq 0}$-gradable $\sigma$-twisted
$\overline{SM(1)}$-module  is isomorphic to an irreducible
subquotient of~\mbox{$M(1,\l) \otimes M^{\pm}$}. By using the  structure of
twisted Zhu's algebra $A_{\sigma}(\overline{SM(1)})$ and the methods
developed in~\cite{AdM-2007}, we can also construct logarithmic
$\sigma$-twisted $\overline{SM(1)}$-modules.
\end{remark}

\section[The $N=1$ Ramond module structure of twisted $V_L \otimes F$-modules]{The $\boldsymbol{N=1}$ Ramond module structure\\ of twisted $\boldsymbol{V_L \otimes F}$-modules}
\label{free-fields}

In this section we shall assume that the reader is familiar with
basic results on twisted representations of lattice vertex
superalgebras. Details can be found in \cite{BK,D2,FLM} and \cite{Xu}.

 We shall use the same notation as in
\cite{AdM-str}. Let $L= {\Z}{\a}$ be a rank one lattice with
nondegenerate form given by
$ \la \a , \a \ra = 2m+1 $, where $m \in {\Zp}$. Let $V_L$ be the
corresponding vertex  superalgebra.

For $i \in {\Z}$, we  set
\begin{gather*}
 \ga_i = \frac{i}{2m+1} \a, \qquad \ga^{R}_i = \frac{\a}{2
(2m+1)} + \frac{i}{2m+1} \a .
\end{gather*}

Then  $\sigma_1 = \exp[ \tfrac{\pi i} {2m+1}\a (0)]$ is a
canonical automorphism of order two of $V_L$. The set
$\{ V_{\ga^{R}_i + L}$,   $i=0, \dots, 2m \} $
provides all the irreducible $\sigma_1$-twisted $V_L$-modules.

\begin{remark}It is important to  notice that $\sigma_1$-twisted
$V_L$-module $(V_{\ga^{R}_i + L}, Y_{\ga^{R}_i + L}) $ can be
constructed from untwisted module $(V_{\ga_i + L}, Y_{\ga_i + L})$
as follows (cf.~\cite{Xu}):
\begin{gather*}
 V_{\ga^{R}_i + L}:= V_{\ga_i + L} \qquad \mbox{as vector space;}
\\
Y_{\ga^{R}_i + L}(\cdot,z):=Y_{\ga_i + L}\big(\Delta\big(\tfrac{\a
}{2(2m+1)},z\big)  \cdot,z\big),
\end{gather*}
where
\[
\Delta(h,z):= z^{h(0)} \exp \left( \sum_{n=1} ^{\infty} \frac{h(n)
}{-n}(-z) ^{-n} \right).
\]
\end{remark}

 Let $F$ be the fermionic vertex operator superalgebra with
central charge $1/2$ and $\sigma_2$ its parity map. Let $M$ be the
  $\sigma_2$ twisted $F$-module (cf.~\cite{FRW}).

Then $\sigma= \sigma_1 \otimes \sigma_2$ is the parity  automorphism
of order two of the vertex superalgebra $V_L \otimes F$ and
\[
  V_{\ga^{R}_i + L} \otimes M , \qquad i=0, \dots, 2m
  \]
are  $\sigma$-twisted $V_L \otimes F$-modules.

Let $(W,Y_W)$ be any  $\sigma$-twisted $V_L \otimes F$-module.
From the Jacobi identity for $\sigma$-twisted modules it follows
that the coef\/f\/icients of
\[
Y_W(\tau,z) = \sum_{n \in {\Z}} G(n) z^{-n-\tfrac{3}{2}} \qquad
\mbox{and} \qquad Y_W(\omega,z) = \sum_{n \in {\Z}} L(n) z^{-n-2}
\]
def\/ine a representation of the $N=1$ Ramond algebra.

Recall that
\[
h^{2i+2,2n+1}= \frac{ ( (2i+2) - (2n+1) (2m+1) ) ^2 -
4m^2 }{8 (2m +1)} + \frac{1}{16}.
\]

\begin{proposition}
Assume that $n \in {\Z}$. Then $e^{{\ga}^{R}_i - n \alpha} \otimes
{\bf 1} ^{\R}$ is a singular vector for the $N=1$ Ramond algebra
$\R$ and
\begin{gather*}
U(\R) \big(e^{{\ga}^{R}_i - n \alpha} \otimes
{\bf 1} ^{\R}\big) \cong L^{\R} \big(c_{2m+1,1}, h^{2i+2,2n+1}\big)\\
\hphantom{U(\R) \big(e^{{\ga}^{R}_i - n \alpha} \otimes {\bf 1} ^{\R}\big)}{}  = L^{\R}
\big(c_{2m+1,1}, h^{2i+2,2n+1}\big)^+ \oplus L^{\R} (c_{2m+1,1},
h^{2i+2,2n+1}) ^- .
\end{gather*}
Moreover,
\[
U(\R) \big(e^{{\ga}^{R}_i - n \alpha} \otimes
{\bf 1} ^{\pm}\big) \cong  L^{\R} \big(c_{2m+1,1}, h^{2i+2,2n+1}\big)^{\pm} .
\]
\end{proposition}

As in \cite{AdM-str} (see also  \cite{IK3} and \cite{MR}) we have
the following result.
\begin{lemma} The (screening)  operator
\[
Q = \mbox{\rm Res}_x Y_W\big(e^{\a}\otimes \phi\big({-}\tfrac{1}{2}\big),x\big)
\]
commutes with the action of $\R$.
\end{lemma}

\begin{remark}
By using generalized  (lattice) vertex algebras and their twisted
representations one can also def\/ine the second screening operator
$\widetilde{Q}$ acting between certain $\sigma$-twisted $V_L
\otimes F$-modules such that
\[
[\widetilde{Q}, {\R}] = 0, \qquad [Q, \widetilde{Q} ]  = 0
\] (for
details and some applications see \cite{AdM-lattice}).
\end{remark}

We shall f\/irst present results on the structure of
$\sigma$-twisted $V_L \otimes F$-modules, viewed as modules for the
$N=1$ Ramond  algebra. Each $V_{L+{\ga}^{R}_i}\otimes M$ is a
direct sum of super Feigin--Fuchs modules via
\[
V_{L+{\ga}^{R}_i} \otimes M=\bigoplus_{n \in \mathbb{Z}} \big(M(1) \otimes
e^{{\ga}^{R}_i +n \alpha}\big) \otimes M.
\]

Since operators $Q ^{j}$, $j \in {\Zp}$, commute with the action
of the Ramond algebra, they  are  actually (Lie superalgebra) intertwiners between
super Feigin--Fuchs modules inside  $V_{L+\ga^{R}_i} \otimes M$.

To simplify the notation, we shall identify $e^{\beta}$ with $e^{\beta} \otimes {\bf 1}^{\R}$ for
every $\beta \in L+ \ga^{R}_i$.

Assume that $0 \le i \le m-1$. If $Q^{j} e^{\ga^{R}_i-n \alpha}$
is nontrivial, it is a singular vector of weight
\[
\mbox{wt}\big(Q^{j} e^{\ga^R_i -n \alpha}\big) = \mbox{wt} \big(e^{\ga^R_i -n \alpha}\big) =h ^{2i+2,2n+1}.
\]

Since $\mbox{wt} \big( e^{ \ga^{R}_i + (j-n) \alpha} \big) > \mbox{wt} \big(
e^{\ga^{R}_i -n \alpha} \big)$ if $j  > 2n   $, we conclude that
\begin{gather*} 
 Q^ {j} e^{\ga^{R}_i-n \alpha} = 0 \qquad \mbox{for} \quad j > 2 n.
\end{gather*} One can similarly see that for $ m \le i \le 2m$:
\begin{gather*}
 Q^ {j} e^{\ga^{R}_i-n \alpha} = 0 \qquad \mbox{for} \quad j > 2 n
+1.
\end{gather*}

The following lemma is useful for constructing singular vectors in
$V_{L+\ga^{R}_i} \otimes M$:
\begin{lemma} \label{non-triv-Ramond} \quad {}
\begin{enumerate}\itemsep=0pt
\item[$(1)$]  $ Q^{2n} e^{\ga^R_i - n\a} \ne 0$ \  for  \  $0 \le i
\le m$.

\item[$(2)$]$ Q^{2n+1} e^{\ga^R_i - n\a} \ne 0$ \ for \  $m+1 \le i
\le 2m$.
\end{enumerate}
\end{lemma}

\begin{proof} The proof uses the results on fusion rules from
Proposition~\ref{fusion-R} and is completely ana\-logous to that of
Lemma~6.1 in~\cite{AdM-str}.
\end{proof}

As in the Virasoro algebra case the $N=1$ Feigin--Fuchs modules are
classif\/ied according to their embedding structure.
For the  purposes of our paper we shall focus only on modules of
certain types (Type 4 and 5 in \cite{IK3}). These modules are either
semisimple (Type 5) or they become semisimple after quotienting with
the maximal semisimple submodule (Type 4).
%

The following result follows directly from Lemma
\ref{non-triv-Ramond} and the structure theory of super Feigin--Fuchs
modules, after some minor adjustments of parameters (cf. Type 4
embedding structure in \cite{IK3}).
\begin{theorem} \label{str-fock-R-1}
Assume that $i \in \{0, \dots, m-1\}$.

\begin{enumerate}\itemsep=0pt
 \item[$(i)$] As an $\R$-module, $V_{L+\ga^{R}_i} \otimes M$  is generated by the family of
singular and cosingular vectors $ \widetilde{\rm Sing}_{i}  \bigcup
\widetilde{\rm CSing}_{i}$, where
\begin{gather*}
  \widetilde{\rm Sing}_{i} =  \big\{ u_i ^{(j, n)} \ \vert \ j, n \in {\N}, \ 0 \le j \le 2n \big\}; \\
  \widetilde{\rm CSing}_{i} =  \big\{ w_i ^ {(j, n)} \ \vert \ n \in {\Zp}, 0 \le j \le 2n-1 \big\}.
  \end{gather*}
These vectors satisfy the following relations:
\[  u_i ^{(j, n)}= Q ^{j} e ^{\ga^{R}_i -n \alpha} , \qquad Q ^{j}
w_i ^{(j, n)} = e ^{\ga^{R}_i + n \a}.
\]
 The submodule generated by singular vectors
$\widetilde{\rm Sing}_i$, denoted by $R\Lambda(i+1)$, is isomorphic to
\[
\bigoplus_{n =0 } ^{\infty} (2n +1) L^{\R}\big(c_{2m+1,1}, h^{2i+2,
2n+1}\big).
\]

\item[$(ii)$] Let $R\Lambda(i+1) ^{\pm}$ denote the submodule of
$V_{L+\ga^{R}_i} \otimes M^{\pm}$ generated by singular vectors
\[
Q ^{j} \big(e ^{\ga^{R}_i -n \alpha} \otimes {\bf 1} ^{\pm}\big), \qquad n \in {\N}, \qquad 0
\le j \le 2n.
\]
Then
\[
 R\Lambda(i+1) ^{\pm } =\bigoplus_{n =0 } ^{\infty} (2n +1) L^{\R}\big(c_{2m+1,1}, h^{2i+2,
2n+1}\big) ^{\pm},
\]
and
\[
 R\Lambda(i+1) = R\Lambda(i+1) ^{+} \oplus
R\Lambda(i+1)^{-}.
\]

\item[$(iii)$] For the quotient module we have
\[
 R \Pi(m-i):=(V_{L+\ga^{R}_i} \otimes M) / R\Lambda(i+1) \cong
  \bigoplus_{n =1 } ^{\infty}
(2 n ) L^{\R}\big(c_{2m+1,1}, h^{2i+2, -2n +1}\big).
\]

Moreover, we have
\[
 R\Pi(m-i) = R\Pi(m-i) ^{+} \oplus R\Pi(m-i)^{-},
 \]
where
\[
 R\Pi(m-i) ^{\pm } = (V_{L+\ga^{R}_i} \otimes M) / R\Lambda(i+1)^{\pm} = \bigoplus_{n =1 } ^{\infty}
(2 n ) L^{\R}\big(c_{2m+1,1}, h^{2i+2, -2n +1}\big) ^{\pm}.
\]
\end{enumerate}
\end{theorem}

\begin{theorem} \label{str-fock-R-2}
Assume that $i \in \{0, \dots, m-1 \}$.

\begin{enumerate}\itemsep=0pt
 \item[$(i)$] As an  $\R$-module, $V_{L+\ga^{R}_{m+i}} \otimes M$  is generated by the family of
singular and cosingular vectors $ \widetilde{\rm Sing}'_{i}  \bigcup
\widetilde{\rm CSing}'_{i}$, where
\begin{gather*}
  \widetilde{\rm Sing}'_{i} = \big\{ u_i^ {'(j, n)} \ \vert \ n \in {\Zp}, 0 \le j \le 2n-1
\big\}; \\
 \widetilde{\rm CSing}'_{i} = \big\{ w_i^{'(j, n)} \ \vert \ j, n \in
{\N}, \ 0 \le j \le 2n \big\}.
\end{gather*}
These vectors satisfy the following relations:
\[
 u_i ^{'(j, n)}= Q ^{j} e ^{\ga^{R}_{m+i} -n \alpha} , \qquad Q
^{j} w_i ^{'(j, n)} = e ^{\ga^{R}_{m+i} + n \a}.
\]
The submodule generated by singular vectors $\widetilde{\rm Sing}_i$
is isomorphic to
\[
 R\Pi(i+1)\cong \bigoplus_{n =1 } ^{\infty} (2 n ) L^{\R}\big(c_{2m+1,1}, h^{2m -2i,
-2n +1}\big).
\]

\item[(ii)] For the quotient module we have
\[
R \Lambda(m -i)\cong (V_{L+\ga^{R}_i} \otimes M) / R \Pi(i+1) \cong
  \bigoplus_{n =0 } ^{\infty}
(2 n+1 )  L^{\R}\big(c_{2m+1,1}, h^{2m -2i, 2n+1}\big).
\]
\end{enumerate}
\end{theorem}

\newpage

\begin{theorem} \label{str-fock-R-3}\quad
\begin{enumerate}\itemsep=0pt
\item[$(i)$] As an  $\R$-module, $V_{L +
\ga^{R}_{2m}} \otimes M$ is completely reducible and  generated by
the family of singular vectors
\[
  \widetilde{\rm Sing}_{2 m} =  \big\{ u_{2 m} ^{(j, n)} : = Q ^{j} e ^{\ga^{R}_{2m} -n \alpha}  \ \vert \ n \in {\Zp}, \  j \in {\N}, \ 0 \le j \le 2n-1 \big\};
\]
 and it is isomorphic to
\[
R\Pi(m+1):=
 V_{L + \ga^{R}_{2m}} \otimes M \cong  \bigoplus_{n=1} ^{\infty} (2n)L^{\R}\big(c_{2m+1,1}, h^{4m+2, -2n
+1}\big).
\]

\item[$(ii)$] Let $R\Pi(m+1) ^{\pm}$ be the submodule of $V_{L +
\ga^{R}_{2m}} \otimes M^{\pm}$ generated by the singular vectors
\[
 Q ^{j} \big(e ^{\ga^{R}_{2m} -n \alpha} \otimes {\bf 1}^{\pm} \big),  \qquad n \in {\Zp}, \qquad  j \in {\N}, \qquad 0 \le j \le
2n-1.
\]
Then
\[
 R\Pi(m+1) ^{\pm } = \bigoplus_{n=1} ^{\infty} (2n)L^{\R}\big(c_{2m+1,1}, h^{4m+2, -2n
+1}\big)^{\pm},
\]
and
\[
 R\Pi(m+1) = R\Pi(m+1) ^{+} \oplus R\Pi(m+1)^{-}.
 \]
\end{enumerate}
\end{theorem}



\begin{remark}
In this section we actually constructed explicitly all the
non-trivial intertwining operators from Proposition \ref{fusion-R}.
\end{remark}

\section[The $\sigma$-twisted $\mathcal{SW}(m)$-modules]{The $\boldsymbol{\sigma}$-twisted $\boldsymbol{\mathcal{SW}(m)}$-modules}
\label{classification}

Since   $\striplet \subset V_L \otimes F $ is $\sigma$-invariant,
then every $\sigma$-twisted $V_L \otimes F$-module is also
a $\sigma$-twisted module for the vertex operator superalgebra
$\striplet$. In this section we shall consider $\sigma$-twisted
$V_L \otimes F$-modules from Section \ref{free-fields} as
$\sigma$-twisted $\striplet$-modules. In what follows we shall
classify all the irreducible $\sigma$-twisted $\striplet$-modules
by using Zhu's algebra~$A_{\sigma}(\striplet)$.

Following  \cite{AdM-triplet} and \cite{AdM-str}, we f\/irst  notice
the following important fact:
\begin{gather}
 Q ^2 e^{-2 \alpha} \in O(\striplet) .
\label{zhu-Q2}
\end{gather}

\begin{proposition} \label{identitet}
Let $v_{\l}^{\pm}$ be the highest weight vector in   $M(1,\l)
\otimes M^{\pm}$. We have
\[
 o\big(Q ^2 e^{-2 \alpha}\big)  v_{\l}^{\pm} =    F_m (t) \,  v_{\l} ^{\pm} ,
 \]
where $t = \la \l, \a \ra$, and
\[
F_m (t) =A_m  { t + m + 1/2 \choose 3 m + 1} { t - 1/2 \choose 3m+1 }, \qquad \mbox{where} \quad A_m = (-1)^m \frac{{2m \choose m} }{{4m+1 \choose m}}.
\]
\end{proposition}

\begin{proof} First we notice that
\begin{gather*}
  Q^{2} e^{-2\a} = w_1 + w_2,
  \end{gather*}
  where
\begin{gather*}
 w_1= \sum_{i=0} ^{\infty } e^{\a}_{-i-1} e^{\a} _i
e^{-2\a}, \qquad
 w_2 =\sum_{i,j=0} ^{\infty} e^{\a}_{i} e^{\a}_j e^{-2\a}
\otimes \phi\big({-}i-\tfrac{1}{2}\big) \phi\big({-}j -\tfrac{1}{2}\big) {\bf 1}.
\end{gather*}
The proof of Proposition~8.3 from \cite{AdM-str} gives that
\begin{gather*}
o(w_1) v_{\l}^{\pm} = {\rm Res}_{x_1} {\rm
Res}_{x_2}(x_2-x_1)^{2m} (1+x_1)^t (1+x_2)^t (x_1 x_2)^{-4m-2}
v_{\l} \\
\phantom{o(w_1) v_{\l}^{\pm}}{} = A_m { t + m   \choose 3 m + 1} { t \choose
3m+1 } v_{\l}^{\pm},
\end{gather*}
 so it remains to examine
$o(w_2)$. Recall Lemma~\ref{cmn}, so that
\[
(x_2 -x_1)G(x_1,x_2)=\tfrac{1}{2} \left( (1+x_1)^{1/2}(1+x_2)^{-1/2}+(1+x_1)^{-1/2}(1+x_2)^{1/2}-2 \right).
\]
Now, we compute
\begin{gather}
  o(w_2) \cdot v_{\l}^{\pm}
=\sum_{i,j \geq 0} o \left( e^{\alpha}_i e^{\alpha}_j e^{-2 \alpha} \otimes \phi\big({-}i-\tfrac{1}{2}\big)\phi\big({-}j-\tfrac{1}{2}\big){\bf 1} \right) \cdot v_{\l} ^{\pm} \nonumber \\
\phantom{o(w_2) \cdot v_{\l}^{\pm}}{}  =\sum_{i,j \geq 0} {\rm Res}_{x_1}{\rm Res}_{x_2} x_1^i x_2^j o\big(Y(e^{\alpha},x_1)Y(e^{\alpha},x_2)e^{-2 \alpha} \big) o\big(\psi\big({-}i-\tfrac{1}{2}\big)\psi\big({-}j-\tfrac{1}{2}\big)\big) \cdot v_{\l} ^{\pm}
 \nonumber \\
\phantom{o(w_2) \cdot v_{\l}^{\pm}}{}
 ={\rm Res}_{x_1} {\rm Res}_{x_2}  \Bigg( \sum_{i,j \geq 0} x_1^i x_2^j  c_{i,j}
{\rm Res}_{x_1} {\rm Res}_{x_2} (x_2-x_1)^{2m+1}(1+x_1)^t \nonumber\\
\phantom{o(w_2) \cdot v_{\l}^{\pm}{\rm Res}_{x_1} {\rm Res}_{x_2}  \Bigg(}{}
\times(1+x_2)^t (x_1x_2)^{-4m-2} v_{\l} ^{\pm} \Bigg) \nonumber \\
\phantom{o(w_2) \cdot v_{\l}^{\pm}}{} ={\rm Res}_{x_1} {\rm Res}_{x_2}(G(x_1,x_2) (x_2-x_1)^{2m+1}(1+x_1)^t (1+x_2)^t )(x_1x_2)^{-4m-2} v_{\l} ^{\pm} \nonumber \\
 \label{-2}
 \phantom{o(w_2) \cdot v_{\l}^{\pm}}{} =-{\rm Res}_{x_1} {\rm Res}_{x_2}(x_2-x_1)^{2m} (1+x_1)^t (1+x_2)^t (x_1 x_2)^{-4m-2} v_{\l} ^{\pm} \\
 \label{-1}
 \phantom{o(w_2) \cdot v_{\l}^{\pm}=}{}+\tfrac{1}{2}{\rm Res}_{x_1} {\rm Res}_{x_2}\left((x_2-x_1)^{2m} (x_1x_2)^{-4m-2} (1+x_1)^{t+1/2} (1+x_2)^{t-1/2} \right) v_{\l} ^{\pm} \!\!\!  \\
 \label{last}
 \phantom{o(w_2) \cdot v_{\l}^{\pm}=}{}
 +\tfrac{1}{2}{\rm Res}_{x_1} {\rm
Res}_{x_2}\left((x_2-x_1)^{2m} (x_1x_2)^{-4m-2} (1+x_1)^{t-1/2}
(1+x_2)^{t+1/2} \right) v_{\l} ^{\pm}.\!\!\!
 \end{gather}
  Now, observe that the
expression in (\ref{-2}) is precisely $-o(w_1) \cdot v_{\l}
^{\pm}$, while (\ref{-1}) and (\ref{last}) are equal.
Consequently,
\begin{gather*}
o(w_1+w_2) \cdot v_{\l} ^{\pm}={\rm Res}_{x_1} {\rm Res}_{x_2}\left((x_2-x_1)^{2m} (1+x_1)^{t-1/2} (1+x_2)^{t+1/2} (x_1x_2)^{-4m-2} \right) v_{\l} ^{\pm}.
\end{gather*}
Now, we expand the generalized rational function in the last formula and obtain
\[
o\big(Q^2e^{-2 \alpha}\big) v_{\l}^{\pm}=\sum_{k=0} ^{2 m } (-1) ^k  { 2 m \choose k } { t +1/2 \choose  4 m + 1 -
k} {t -1/2 \choose  2 m + 1+k} v_{\l}^{\pm}.
\] The sum in the last
formula can be evaluated as in \cite{AdM-triplet} and
\cite{AdM-str}. We have
\[
\sum_{k=0} ^{2 m } (-1) ^k  { 2 m \choose k } { t +1/2 \choose  4 m + 1 -
k} {t -1/2 \choose  2 m + 1+k} =A_m  { t + m + 1/2 \choose 3 m + 1} { t - 1/2 \choose 3m+1 },
\]
where $A_m$ is above. The proof follows.
\end{proof}

A direct consequence of Proposition~\ref{identitet} and relation
(\ref{zhu-Q2}) is the following important result:

\begin{theorem} \label{vazna-rel-ramond}
  In Zhu's  algebra $A_{\sigma}(\striplet)$ we have the following relation
\[
f^{\R}_m([\omega]) = 0,
\]
 where
\begin{gather*}
  f^{\R}_m(x) =  \prod_{i=0}^{3 m} \big(x - h^{2i+2,1}\big) =  \left( \prod_{i=0}^{m-1} \big(x-h^{2i+2,1}\big)^{2} \right) \left(\prod_{i=2m} ^{3 m}\big(x-h^{2i+2,1}\big)\right) .
\end{gather*}
\end{theorem}

In parallel with \cite{AdM-str} we conjecture that $f^{\R}_m(x)$ is in fact the minimal polynomial of $[\omega]$ in $A_{\sigma}(\striplet)$.

We have the following irreducibility result. The proof is similar
to that  of Theorem 3.7 in \cite{AdM-triplet}.

\begin{theorem} \label{ireducibilni-moduli-ramond-nongraded} \qquad{}
\begin{enumerate}\itemsep=0pt

\item[$(1)$] For every $0 \leq i \leq m-1$,  $R\Lambda(i+1) ^{\pm}$
are $\N$-gradable irreducible  $\sigma$-twisted
$\mathcal{SW}(m)$-modules and   the top components $R\Lambda(i+1)
^{\pm} (0)$ are $1$-dimensional irreducible
  $A_{\sigma}(\mathcal{SW}(m))$-modules.

  \item[$(1')$] For every $0 \leq i \leq m-1$,  $R\Lambda(i+1)$ is a graded irreducible  $\sigma$-twisted $\mathcal{SW}(m)$-module,
and  its  top component $R\Lambda(i+1) (0)$ is a graded irreducible
  $A_{\sigma}(\mathcal{SW}(m))$-module.

\item[$(2)$]  For every $0 \leq j \leq m$ $,  R\Pi(j+1) ^{\pm}$ are
irreducible   $\N$-gradable $\sigma$-twisted
$\mathcal{SW}(m)$-modules and the top components  $S\Pi(j+1)
^{\pm} (0)$ are irreducible  $2$-dimensional
$A_{\sigma}(\mathcal{SW}(m))$-modules.

\item[$(2')$]  For every $0 \leq j \leq m$ $,  R\Pi(j+1)$ is a
  graded irreducible $\sigma$-twisted
$\mathcal{SW}(m)$-module and its  top component $S\Pi(j+1) (0)$ is
a graded irreducible  $A_{\sigma}(\mathcal{SW}(m))$-module.
\end{enumerate}
\end{theorem}

\begin{corollary} The minimal polynomial of $[\omega]$ is divisible by
  \[
  \prod_{i=0}^{m-1} (x-h^{2i+2,1}) \left(\prod_{i=2m} ^{3 m}(x-h^{2i+2,1})\right).
  \]
Moreover, both $[\omega]$ and $[\tau]$ are units in
$A_{\sigma}(\striplet)$. Consequently, $\striplet$ has no
supersymmetric sector (i.e., there is no $\sigma$-twisted
$\striplet$-modules of highest weight $\frac{c_{2m+1,1}}{24}$).
\end{corollary}

By using similar arguments as in  \cite{AdM-triplet} and \cite{AdM-str} we have the following
result on the structure of Zhu's algebra $A_{\sigma}(\striplet)$ (as in the
proof of Theorem 4.1 we see that $[\tau]$ is central).
\begin{proposition} \label{zhu-komut-r}
The associative algebra $A_{\sigma} (\striplet)$ is  generated
by $[E]$, $[H]$, $[F]$, $[\tau]$, $[\widehat{E}]$,
$[\widehat{H}]$, $[\widehat{F}]$ and $[\omega]$.  The following
relations hold:
\begin{gather*}
 [\tau] \quad \mbox{and} \quad  [\omega]  \quad
\mbox{are in the center of} \ A_{\sigma}(\striplet),   \\
 [\tau] ^2 = [\omega] - \frac{c_{2m+1,1}}{24}, 
 \\
 [\tau]*[X]=\frac{-\sqrt{2m+1}}{2} [\widehat{X}], \qquad \mbox{for}
\quad X \in \{E,
F, H\}, \\
[\widehat{H}]*[\widehat{F}]-[\widehat{F}]*[\widehat{H}]=-2q([\omega])[\widehat{F}], \\
 [\widehat{H}]*[\widehat{E}]-[\widehat{E}]*[\widehat{H}]=2q([\omega])[\widehat{E}], \\
[\widehat{E}]*[\widehat{F}]-[\widehat{F}]*[\widehat{E}]=-2q([\omega])[\widehat{H}],
\end{gather*}
 where $q$ is a certain polynomial.
 \end{proposition}

Equipped with all these results we are now ready to classify irreducible $\sigma$-twisted $\striplet$-modules.

\begin{remark}
By using same arguments as in Proposition 5.6 of~\cite{ABD} we have that
for $C_2$-cof\/inite SVOAs, every week (twisted) module is admissible (see also~\cite{DZ}). Thus, the classif\/ication of
irreducible $\sigma$-twisted $\striplet$-modules reduces to classif\/ication of irreducible $\mathbb{Z}_{\geq 0}$-gradable modules.
\end{remark}

\begin{theorem} \label{class-ired-zhu}\qquad{}
\begin{enumerate}\itemsep=0pt
\item[$(i)$]The set
\[
\{ R\Pi(i)^{\pm}(0),  : 1 \leq i \leq m+1  \}
\cup \{ R\Lambda(i)^{\pm}(0) : 1 \leq i \leq m\}
\] provides, up to
isomorphism, all irreducible modules for  Zhu's algebra
$A_{\sigma}(\mathcal{SW}(m))$.

\item[$(ii)$] The set
\[
\{ R\Pi(i)(0),  : 1 \leq i \leq m+1 \} \cup
\{ R\Lambda(i)(0) : 1 \leq i \leq m\}
\] provides, up to
isomorphism, all ${\Z}_2$-graded irreducible modules for  Zhu's
algebra \linebreak $A_{\sigma}(\mathcal{SW}(m))$.
\end{enumerate}
\end{theorem}

\begin{proof}
The proof  is similar to those of Theorem 3.11 in
\cite{AdM-triplet} and Theorem 10.3 in \cite{AdM-str}.
 Assume that $U$ is an irreducible $A_{\sigma}(\mathcal{SW}(m))$-module. Relation
$f_m ^{\R} ([\omega]) = 0$ in $A_{\sigma}(\mathcal{SW}(m))$
implies that
\[
L(0) \vert U = h^{2i+2,1} \, \mbox{Id}, \qquad \mbox{for} \quad i \in \{0, \dots, m-1 \} \cup \{2m, \dots, 3m\}.
\]
Moreover, since $[\tau]$ is   in the center of
$A_{\sigma}(\mathcal{SW}(m))$, we conclude that $G(0)$ also acts on
$U$ as a~scalar. Then relation $[\tau] ^2 = [\omega] - c_{2m+1,1} /
24$ implies that
\[
  G(0) \vert U = \frac{1/2 + i -m}{\sqrt{2(2m+1)}}  \, \mbox{Id} \qquad
\mbox{or} \qquad G(0) \vert U = -\frac{1/2 + i -m}{\sqrt{2(2m+1)}}
\, \mbox{Id}.
\]
Assume f\/irst that $i= 2m + j$ for $ 0 \le j \le m$. By combining
Propositions \ref{zhu-komut-r} and Theorem
\ref{ireducibilni-moduli-ramond-nongraded} we have that
$q(h^{2i+2,1}) \ne 0$. Def\/ine
\[
 e= \frac{ 1}{ \sqrt{2} q(h^{2 i+2,1})} [\widehat{E}], \qquad f= -\frac{1}{ \sqrt{2}q(h^{2i+2,1})} [\widehat{F}], \qquad h= \frac{1}{q(h^ {2i+2,1})} [\widehat{H}] .
 \]

Therefore $U$ carries the structure of an irreducible,
$\goth{sl}_2$-module with the property that $e^2 = f^2 = 0$  and
$h \ne 0$ on $U$. This easily implies that $U$ is a
$2$-dimensional irreducible $\goth{sl}_2$-module. Moreover,  as
an $A_{\sigma}(\mathcal{SW}(m))$-module $U$ is isomorphic to
either  $R\Pi(m+1-j)^{+}(0)$ or $R\Pi(m+1-j)^{-}(0)$.

In the case  $ 0 \le i \le m-1$, as in \cite{AdM-triplet} we prove
  that $U \cong R\Lambda(i+1)^{+} (0)$ or $U \cong R\Lambda(i+1)^{-} (0)$.

  Let us now prove the second assertion. Let $ N = N^{0} \oplus
  N^{1}$ be the graded irreducible
  $A_{\sigma}(\striplet)$-module. As above, we have that
\[
L(0) \vert N= h^{2i+2,1} \, \mbox{Id}, \qquad \mbox{for} \quad i \in \{0, \dots, m-1 \} \cup \{2m, \dots, 3m\}.
\]

  Then $N= N^{+} \oplus N^{-}$, where
\begin{gather*}
 N^{\pm} = \mbox{span}_{\C} \left\{ v \pm \frac{\sqrt{2(2m+1)}}{1/2 + i -m} G(0) v  \ \vert v \in
  N^{0} \right\} \qquad \!\!\mbox{and}\!\!\qquad
 G(0) \vert N^{\pm} = \pm \frac{1/2 + i -m}{\sqrt{2(2m+1)}}  \, \mbox{Id} .
\end{gather*}
By using assertion $(i)$, we easily get that $N^{\pm} =
R\Lambda(i+1)^{\pm} (0)$ (if $ 0 \le i \le m-1$)  and  $N ^{\pm} =
R\Pi(3m+1-i)^{\pm}(0) $ (if $ 2m \le i \le 3m$).
  \end{proof}

\begin{theorem} \label{class-ired}\qquad
\begin{enumerate}\itemsep=0pt
\item[$(i)$] The set
\[
\{ R\Pi(i)^{\pm}   : 1 \leq i \leq m+1  \}
\cup \{ R\Lambda(i)^{\pm} : 1 \leq i \leq m\}
\]
 provides, up to
isomorphism, all irreducible $\sigma$-twisted
$\striplet$-modules.

\item[$(ii)$] The set
\[
\{ R\Pi(i)   : 1 \leq i \leq m+1  \} \cup \{
R\Lambda(i) : 1 \leq i \leq m\}
\] provides, up to isomorphism, all
${\Z}_2$-graded irreducible $\sigma$-twisted
$\striplet$-modules.
\end{enumerate}
\end{theorem}

So the vertex operator algebra $\striplet$ contains only f\/initely
many irreducible modules. But one can easily see that modules
$V_{L+\ga^{R}_{i}} \otimes M$ and $V_{L+\ga^{R}_{m+i}} \otimes M$
($0 \le i \le m-1$) constructed in Theorems \ref{str-fock-R-1} and
\ref{str-fock-R-2} are not completely reducible. Thus we have:

\begin{corollary} The vertex operator superalgebra $\striplet$ is
not $\sigma$-rational, i.e., the category of $\sigma$-twisted
$\striplet$-modules is not semisimple.
\end{corollary}

\begin{remark}
In our    forthcoming paper \cite{AdM-lattice}   we shall prove that
$\striplet$ also contains logarithmic $\sigma$-twisted
representations.
\end{remark}

\section[Modular properties of characters of $\sigma$-twisted  $\mathcal{SW}(m)$-modules]{Modular properties of characters\\ of $\boldsymbol{\sigma}$-twisted  $\boldsymbol{\mathcal{SW}(m)}$-modules}
 \label{characters}

We f\/irst introduce some basic modular forms
needed for description of irreducible twisted $\striplet$
characters. The Dedekind $\eta$-function is usually def\/ined as the
inf\/inite product
\[
\eta(\tau)=q^{1/24} \prod_{n=1}^\infty (1-q^n),
\]
an automorphic form of weight $\frac{1}{2}$. As usual in all these
formulas $q=e^{2 \pi i \tau}$, $\tau \in \mathbb{H}$. We also
introduce
\begin{gather*} 
\goth{f}(\tau) = q^{-1/48}
\prod_{n=0
}^\infty (1+q^{n+1/2}), \\
\goth{f}_1(\tau) = q^{-1/48} \prod_{n=1}^\infty (1-q^{n-1/2}),
\\
\goth{f}_2(\tau) = q^{1/24} \prod_{n=1}^\infty (1+q^n).
\end{gather*}

Let us recall Jacobi $\Theta$-function
\[
\Theta_{j,k}(\tau,z)=\sum_{n \in \mathbb{Z}} q^{k(n+\frac{j}{2k})^2} e^{2 \pi i k z (n+ \frac{j}{2k} ) },
\]
where $j,k \in \frac{1}{2}\mathbb{Z}$.
If $j \in \mathbb{Z}+\frac{1}{2}$ or $k \in \mathbb{N}+\frac{1}{2}$ it will be useful to use the formula
\begin{gather*} 
\Theta_{j,k}(\tau,z)=\Theta_{2j,4k}(\tau,z)+\Theta_{2j-4k,4k}(\tau,z).
\end{gather*}
Observe that
\[
\Theta_{j+2k,k}(\tau,z)=\Theta_{j,k}(\tau,z).
\]
We also let
\[
\partial \Theta_{j,k}(\tau,z):=\frac{1}{\pi i} \frac{d}{dz} \Theta_{j,k}(\tau,z)=\sum_{n \in \mathbb{Z}}
(2kn+j)q^{(2kn+j)^2/4k}.
\] Related $\Theta$-functions needed for
supercharacters are
\[
G_{j,k}(\tau,z)=\sum_{n \in \mathbb{Z}} (-1)^n q^{k(n+\frac{j}{2k})^2} e^{2 \pi i k z (n+ \frac{j}{2k} ) },
\]
where $j \in \mathbb{Z}$ and $k \in \frac{1}{2} \mathbb{N}$. It is easy to see that
\begin{gather*} 
G_{j,k}(\tau,z)=\Theta_{2j,4k}(\tau,z)-\Theta_{2j-4k,4k}(\tau,z).
\end{gather*}
Eventually, we shall let $z=0$ so before we introduce
\[
\Theta_{j,k}(\tau):=\Theta_{j,k}(\tau,0), \qquad G_{j,k}(\tau):=G_{j,k}(\tau,0).
\]
Similarly, we def\/ine $\partial \Theta_{j,k}(\tau)$ and $\partial G_{j,k}(\tau)$.

The following formulas will be useful:
\begin{gather*} 
  \Theta_{j,k}\left(\frac{-1}{\tau},\frac{z}{\tau}\right)=\frac{\sqrt{-i \tau}}{\sqrt{2k}} e^{\frac{4 \pi i k z^2}{\tau}}\left(\sum_{j'=0}^{4k-1} e^{-\pi i j' j/k} \Theta_{2j',4k}(\tau,z) \right), \\
  \partial \Theta'_{j,k}\left(\frac{-1}{\tau},\frac{z}{\tau}\right)=\frac{\sqrt{-i \tau}}{\sqrt{2k}} e^{\frac{4 \pi i k z^2}{\tau}} {8 k z} \left(\sum_{j'=0}^{4k-1} e^{-\pi i j' j/k} \Theta_{2j',4k}(\tau,z) \right)\nonumber  \\
\hphantom{\partial \Theta'_{j,k}\left(\frac{-1}{\tau},\frac{z}{\tau}\right)=}{}  + \tau \frac{\sqrt{-i \tau}}{\sqrt{2k}} e^{\frac{4 \pi i k z^2}{\tau}} \left(\sum_{j'=0}^{4k-1} e^{-\pi i j' j/k} \partial \Theta_{2j',4k}(\tau,z) \right).
\end{gather*}
If we let now $z=0$, we obtain
\begin{gather} \label{theta-1}
\Theta_{j,k}\left(\frac{-1}{\tau}\right)=\frac{\sqrt{-i \tau}}{\sqrt{2k}}\sum_{j'=0}^{4k-1} e^{\frac{-\pi i j j'}{k}} \Theta_{2j',4k}(\tau,0),
\end{gather}
and
\begin{gather} \label{theta-2}
\partial \Theta'_{j,k}\left(\frac{-1}{\tau}\right)= \tau \frac{\sqrt{-i \tau}}{\sqrt{2k}}\sum_{j'=0}^{4k-1} e^{\frac{-\pi i j j'}{k}} \partial \Theta_{2j',4k}(\tau,0).
\end{gather}

For $j \in \mathbb{Z}$ and $k \in \mathbb{N}+\frac{1}{2}$ we now have:
\begin{gather*} 
 \Theta_{j,k}\left(\frac{-1}{\tau}\right)=\sqrt{\frac{-i
\tau}{2k}}\sum_{j'=0}^{2k-1} e^{-i \pi j j'/k} \Theta_{j',k}(\tau), \\
 \Theta_{j,k}(\tau+1)=e^{i \pi j^2/2k}G_{j,k}(\tau), \\
 (\partial \Theta)_{j,k}(\tau+1)=e^{i \pi j^2/2k}(\partial
G)_{j,k}(\tau), \\
 (\partial \Theta)_{j,k}\left(\frac{-1}{\tau}\right)=\tau \sqrt{-i \tau/2k}
\sum_{j'=1}^{2k-1} e^{-i \pi j j'/k} (\partial
\Theta)_{j',k}(\tau),
\end{gather*}
where we used (\ref{theta-1}) and (\ref{theta-2}), together with  $\Theta_{j,k}(\tau)=\Theta_{j+2k,k}(\tau)$  and
$\partial \Theta_{j,k}(\tau)=\partial \Theta_{j+2k,k}(\tau)$.

For $j \in \mathbb{Z}+\frac{1}{2}$ and $k \in \mathbb{N}+\frac{1}{2}$ the transformation formulas are slightly
dif\/ferent:
\begin{gather*} 
 \Theta_{j,k}\left(\frac{-1}{\tau}\right)=\sqrt{\frac{-i
\tau}{2k}}\sum_{j'=0}^{2k-1} e^{-i \pi j j'/k} G_{j',k}(\tau), \\
 \Theta_{j,k}(\tau+1)=e^{i \pi j^2/2k}\Theta_{j,k}(\tau), \\
 (\partial \Theta)_{j,k}(\tau+1)=e^{i \pi j^2/2k}(\partial
\Theta)_{j,k}(\tau), \\
  (\partial \Theta)_{j,k}\left(\frac{-1}{\tau}\right)=\tau \sqrt{-i \tau/2k}
\sum_{j'=1}^{2k-1} e^{-i \pi j j'/k} (\partial
G)_{j',k}(\tau).
\end{gather*}

For a (twisted) vertex operator algebra module $W$ we def\/ine its
graded-dimension or simply character
\[
\chi_{W}(\tau)={\rm tr}|_{W} q^{L(0)-c/24}.
\]
If $V=L^{\NS}(c_{2m+1,0},0)$ and
$W=L^{\R}(c_{2m+1,0},h^{2i+2,2n+1})$, then (see \cite{IK3}, for
instance)
 \begin{gather} \label{irr-r}
 \chi_{L^{\R}(c_{2m+1,1},h^{2i+2,2n+1})}(\tau)=2 q^{\frac{m^2}{2(2m+1)}-\frac{1}{16}}\frac{\goth{f}_2
 (\tau)}{\eta(\tau)}\left(q^{h^{2i+2,2n+1}}-q^{h^{2i+2,-2n-1}}\right).
\end{gather}

By combining Theorems~\ref{str-fock-R-1}, \ref{str-fock-R-2} and
\ref{str-fock-R-3}, and formula (\ref{irr-r}) we obtain

\begin{proposition} \label{char-irr-mod-r} For $i=0,\dots,m-1$
\begin{gather} \label{sl-char}   \chi_{ R\Lambda(i+1)}(\tau)=
2 \frac{\goth{f}_2(\tau)}{\eta(\tau)}\left(\frac{2i+2}{2m+1}
\Theta_{m-i-\tfrac{1}{2},\frac{2m+1}{2}}(\tau)+\frac{2}{2m+1}(\partial
\Theta)_{m-i-\tfrac{1}{2},\frac{2m+1}{2}}(\tau) \right), \\
 \label{pi-char}
 \chi_{R \Pi(m-i)}(\tau)= 2
\frac{\goth{f}_2(\tau)}{\eta(\tau)}\left(\!\frac{2m-2i-1}{2m+1}
\Theta_{m-i-\tfrac{1}{2},\frac{2m+1}{2}}(\tau)-\frac{2}{2m+1}(\partial
\Theta)_{m-i-\tfrac{1}{2},\frac{2m+1}{2}}(\tau)\!\right)\!. \!\!\!\!\!
\end{gather}
 Also, \begin{gather*}
  \chi_{R \Pi(m+1)}(\tau)=2
\frac{\goth{f}_2(\tau)}{\eta(\tau)} \Theta_{m+
\tfrac{1}{2},\frac{2m+1}{2}}(\tau).
\end{gather*}
\end{proposition}

As in \cite{AdM-str}, the characters of irreducible
$\sigma$-twisted $\striplet$-modules can be described by using
characters of irreducible modules for the triplet vertex algebra
$\triplet$, where $p= 2 m +1$. Let $\Lambda(1), \dots, \Lambda
(p)$, $\Pi(1), \dots , \Pi(p)$ be the irreducible
$\triplet$-module. We have the following result:

\begin{proposition}\qquad{}
\begin{enumerate}\itemsep=0pt
\item[$(i)$] For $0 \le i \le m-1$, we have:
\[
 \chi_{ R\Lambda(i+1)}(\tau) = 2 \frac{\chi_{ \Lambda(2
i+2)}(\tau/2)}{\goth{f}(\tau)}.
\]

\item[$(ii)$] For $0 \le i \le m$, we have:
\[
 \chi_{ R\Pi(m+1 -i)}(\tau) = 2\frac{\chi_{ \Pi(2m - 2i+1)}(\tau/2)}{\goth{f}(\tau)}.
 \]
\end{enumerate}
\end{proposition}

Now, we recall also formulas for irreducible $\striplet$ characters and supercharacters  obtained in \cite{AdM-str}.

For $i=0,\dots,m-1$
\begin{gather} \label{orsl-char}   \chi_{S
\Lambda(i+1)}(\tau)=\frac{\goth{f}(\tau)}{\eta(\tau)}\left(\frac{2i+1}{2m+1}
\Theta_{m-i,\frac{2m+1}{2}}(\tau)+\frac{2}{2m+1}(\partial
\Theta)_{m-i,\frac{2m+1}{2}}(\tau) \right), \\
  \label{orpi-char} \chi_{S
\Pi(m-i)}(\tau)=\frac{\goth{f}(\tau)}{\eta(\tau)}\left(\frac{2m-2i}{2m+1}
\Theta_{m-i,\frac{2m+1}{2}}(\tau)-\frac{2}{2m+1}(\partial
\Theta)_{m-i,\frac{2m+1}{2}}(\tau)\right).
\end{gather}
 Also, \begin{gather}
  \label{ss-char}    \chi_{S
\Lambda(m+1)}(\tau)=\frac{\goth{f}(\tau)}{\eta(\tau)}
\Theta_{0,\frac{2m+1}{2}}(\tau).
\end{gather}

For supercharacters we have: for $i=0,\dots,m-1$
\begin{gather}\label{sl-schar}   \chi^F_{S \Lambda(i+1)}(\tau)=
 \frac{\goth{f}_1(\tau)}{\eta(\tau)} \biggl(
\frac{2i+1}{2m+1}G_{m-i,\frac{2m+1}{2}}(\tau)+ \frac{2}{2m+1} (\partial G)_{m-i,\frac{2m+1}{2}}(\tau) \biggr),\\
 \label{p1-schar}   \chi^F_{S \Pi(m-i)}(\tau)=
  \frac{\goth{f}_1(\tau)}{\eta(\tau)} \biggl(
\frac{2m-2i}{2m+1}G_{m-i,\frac{2m+1}{2}}(\tau)- \frac{2}{2m+1} (\partial G)_{m-i,\frac{2m+1}{2}}(\tau) \biggr).
\end{gather}
Also, \begin{gather}
 \label{ss-schar} \chi_{S
\Lambda(m+1)}^F(\tau)=\frac{\goth{f}_1(\tau)}{\eta(\tau)}
G_{0,\frac{2m+1}{2}}(\tau).
\end{gather}

These characters and supercharacters  can be expressed by using
characters of $W(2m+1)$-modules.

\begin{proposition} \label{SWW}\qquad {}{\samepage
\begin{itemize}\itemsep=0pt
\item[$(i)$] For  $0 \leq i \leq m$, we have
\[
\chi_{S\Lambda(i+1)}(\tau)=\frac{\chi_{\Lambda(2i+1)}(\frac{\tau}{2})}{\goth{f}_2(\tau)}, \qquad
\chi_{S\Lambda(i+1)}
^F(\tau)=\frac{\chi_{\Lambda(2i+1)}(\frac{\tau+1}{2})}{\goth{f}_2(\tau)}.
\]
\item[$(ii)$] For $0 \leq i \leq m-1$, we also have
\[
\chi_{S\Pi(m-i)}(\tau)=\frac{\chi_{\Pi(2m-2i)}(\frac{\tau}{2})}{\goth{f}_2(\tau)}, \qquad
\chi_{S\Pi(m-i)}
^F(\tau)=\frac{\chi_{\Pi(2m-2i)}(\frac{\tau+1}{2})}{\goth{f}_2(\tau)}.
\]
\end{itemize}
Here $\Lambda(i)$ and $\Pi(2m+2-i)$, $i=1, \dots, 2m+1$, are
irreducible $\mathcal{W}(2m+1)$-modules~{\rm \cite{AdM-triplet}}.}
\end{proposition}

By combining   transformation formulas for $\Theta_{j,k}(\tau)$,
$\partial \Theta_{j,k}(\tau)$, formulas
(\ref{sl-char})--(\ref{ss-char}), (\ref{sl-schar})--(\ref{ss-schar})
and Proposition~\ref{char-irr-mod-r} we obtain second main result of
our paper.

\begin{theorem} \label{thm-main-2} The $SL(2,\mathbb{Z})$ closure, called $\mathcal{H}$, of the vector space determined by $\striplet$-characters, $\striplet$-supercharacters and
$\sigma$-twisted $\striplet$-characters has the following basis:
\begin{gather*} 
  \frac{\goth{f}_1(\tau)}{\eta(\tau)} G_{0,\frac{2m+1}{2}}(\tau),  \qquad  \frac{\goth{f}(\tau)}{\eta(\tau)} \Theta_{0,\frac{2m+1}{2}}(\tau), \qquad   \frac{\goth{f}_2(\tau)}{\eta(\tau)}
\Theta_{m+ \tfrac{1}{2},\frac{2m+1}{2}}(\tau), \\
  \frac{\goth{f}_1(\tau)}{\eta(\tau)} G_{m-i,\frac{2m+1}{2}}(\tau),  \qquad \frac{\goth{f}(\tau)}{\eta(\tau)}\Theta_{m-i,\frac{2m+1}{2}}(\tau),  \qquad  \frac{\goth{f}_2(\tau)}{\eta(\tau)} \Theta_{m-i-\tfrac{1}{2},\frac{2m+1}{2}}(\tau) ,  \\
  \frac{\goth{f}_1(\tau)}{\eta(\tau)} \partial G_{m-i,\frac{2m+1}{2}}(\tau),  \qquad \frac{\goth{f}(\tau)}{\eta(\tau)} \partial \Theta_{m-i,\frac{2m+1}{2}}(\tau), \qquad  \frac{\goth{f}_2(\tau)}{\eta(\tau)}  \partial \Theta_{m-i-\tfrac{1}{2},\frac{2m+1}{2}}(\tau) , \\
  \tau \frac{\goth{f}_1(\tau)}{\eta(\tau)} \partial G_{m-i,\frac{2m+1}{2}}(\tau),  \qquad \tau \frac{\goth{f}(\tau)}{\eta(\tau)} \partial \Theta_{m-i,\frac{2m+1}{2}}(\tau),  \qquad \tau \frac{\goth{f}_2(\tau)}{\eta(\tau)}  \partial \Theta_{m-i-\tfrac{1}{2},\frac{2m+1}{2}}(\tau),
\end{gather*}
where $i=0,\dots,m-1$. In particular, the space is $9m+3$ dimensional.
\end{theorem}

\subsection[Modular differential equations for $\sigma$-twisted $\striplet$ characters]{Modular dif\/ferential equations for $\boldsymbol{\sigma}$-twisted $\boldsymbol{\striplet}$ characters}

Let us recall classical $SL(2,\mathbb{Z})$ Eisenstein series ($k \geq 1$):
\[
G_{2k}(\tau)=\frac{-B_{2k}}{(2k)!}+\frac{2}{(2k-1)!}\sum_{n=1}^\infty \frac{n^{2k-1}q^n}{1-q^n},
\]
and certain linear combination of level 2 Eisenstein series ($k \geq 1$):
\begin{gather*}  G_{2k,1}(\tau)=\frac{B_{2k}}{(2k)!}+ \frac{2}{(2k-1)!} \sum_{n \geq 1}
\frac{n^{2k-1}q^n}{1+q^n},   \\
  {G}_{2k,0}(\tau)=\frac{B_{2k}(1/2)}{(2k)!}+\frac{2}{(2k-1)!}\sum_{n=1}^\infty \frac{(n-1/2)^{2k-1}q^{n-1/2}}{1+q^{n-1/2}},
\end{gather*}
where $B_{2k}(x)$ are the Bernoulli polynomials, and $B_{2k}$ are the Bernoulli numbers.

A modular dif\/ferential equation is an (ordinary) dif\/ferential equation of the form:
\[
\left( q \frac{d}{dq} \right)^k y(q)+\sum_{j=0}^{k-1} H_j(q) \left(d \frac{d}{dq} \right)^i y(q)=0,
\]
where  $H_j(q)$ are polynomials in Eisenstein series $G_{2i}$, $i \geq 1$, such that the vector
space of solutions is modular invariant.

It is known (cf.~\cite{Z}) that $C_2$-cof\/initeness condition leads
to certain modular dif\/ferential equation satisf\/ied by irreducible
characters ${\rm tr}_M q^{L(0)-c/24}$. In some instances the
degree of this dif\/ferential equation is bigger than the number of
irreducible characters.  So it is not clear what $k$ should be in
general. In the case of Virasoro minimal models, the degree of the
modular dif\/ferential equation is precisely the number of (linearly
independent) irreducible characters.

If $V$ is a vertex operator superalgebra  one can also get modular dif\/ferential equation satisf\/ied by ordinary characters but with respect to the subgroup $\Gamma_{\theta} \subset SL(2,\mathbb{Z})$, where $H_i$ are polynomials in $G_{2i}$ and $G_{2i,0}$ (see~\cite{M} for the precise  statement in the case of $N=1$ minimal models).

In \cite{AdM-ptraces} we proved that the $C_2$-cof\/initeness for the super triplet  vertex algebra $\striplet$ gives rise to a dif\/ferential equation of order $3m+1$ satisf\/ied by $2m+1$ irreducible characters (additional $m$ solutions can be interpreted as certain
{\em pseudotraces}). By applying arguments similar to those in \cite{AdM-ptraces}  it is not hard to prove

\begin{theorem} \label{diff-equation} The irreducible $\sigma$-twisted $\striplet$ characters satisfy the differential equation of the form
\[
\left(d \frac{d}{dq} \right)^{3m+1} y(q)+\sum_{j=0}^{3m} \tilde{H}_{j,0}(q)\left(q \frac{d}{dq} \right)^i y(q)=0,
\]
where $\tilde{H}_{j,0}(q)$ are certain polynomials in $G_{2i}$ and $G_{2i,1}$.
\end{theorem}

As in the case of $\triplet$-modules and ordinary $\striplet$-modules,  we expect that additional $m$ linearly independent solutions in Theorem \ref{diff-equation} have interpretation
in terms of  $\sigma$-twisted pseudotraces (cf.~Conjecture \ref{conj-twist-modular} below).

\section{Conclusion}

Here we gather a few more-or-less expected conjectures (especially in view of our earlier work~\cite{AdM-triplet} and~\cite{AdM-str}).

The f\/irst one is concerned about the structure of $A_{\sigma}(V)$. As in the case of $A(\triplet)$ and $A(\striplet)$, from our analysis in Chapter 6, we can  show that in fact
\[
{\rm dim}(A_{\sigma}(\striplet)) \geq 10m+8.
\]
But this is well below the conjectural dimension, because $10m+8$-dimensional part cannot control possible  logarithmic modules. Thus, as in the case of ordinary $\striplet$-modules, we expect
to have  $2m$ non-isomorphic (non-graded) logarithmic modules with two-dimensional top component. This then leads to the following conjecture:

\begin{conjecture} \label{conj-twist} For every $m \in \mathbb{N}$,
\[
{\rm dim}(A_{\sigma}(\striplet))=12m+8.
\]
\end{conjecture}

If we assume the existence of $m$ logarithmic modules so that Conjecture \ref{conj-twist} holds true, then the following fact is expected.
\begin{conjecture} \label{conj-twist-modular}
The vector space of (suitably defined) generalized $\sigma$-twisted characters is $3m+1$-dimensional.
\end{conjecture}

Thus,  generalized $\striplet$ characters, supercharacters and $\sigma$-twisted characters together should give rise to  a $9m+3$-dimensional modular invariant space.

\subsection*{Acknowledgments}

The second author was partially supported by NSF grant DMS-0802962.

\pdfbookmark[1]{References}{ref}
\LastPageEnding

\end{document}